\RequirePackage{fix-cm}
\documentclass[smallextended]{svjour3}
\smartqed  % flush right qed marks, e.g. at end of proof
\usepackage{graphicx}
\usepackage{epsfig,amsfonts,amssymb,latexsym,amsmath}

 \journalname{my journal}
\begin{document}
\title{On the Numerical Solution of Non-linear First Order Ordinary Differential Equation Systems}%\thanks{Grants or other notes
\subtitle{ Applications to a Flight Mechanics Problem}

\titlerunning{On the Numerical Solution of ODEs}        % if too long for running head

\author{Fabio Silva Botelho }
       % Second Author %etc.
%}


\institute{Fabio Silva Botelho \at Department of Mathematics \\
              Federal University of Santa Catarina, SC - Brazil \\
              Tel.: +55-48-3721-3663\\
              \email{fabio.botelho@ufsc.br}  }        %  \\
          % \and
           %S. Author \at
              %second address
%}


\maketitle

\begin{abstract}
In this article, firstly we develop a method of solution for a type of difference equations,
 applicable to solve approximately a class of first order  ordinary differential equation systems.

 In a second step, we apply the results obtained to solve a non-linear two point boundary value problem relating a
 flight mechanics model. We highlight the algorithm obtained seems to be  robust and of easy
 computational  implementation.
  \\
%Insert your abstract here. Include keywords, PACS and mathematical
%subject classification numbers as needed.
\keywords{Ordinary differential equations \and Two point boundary value problem \and Flight mechanics} %First keyword \and Second keyword \and More}
\subclass{34L30 \and 39A60}%\and MSC code2 \and more}
\end{abstract}
\section{Introduction}

In this article we develop an algorithm to solve a class of first order non-linear ordinary differential equations.

We start by presenting a general procedure for solving the linearized equations, and in a second step, we
apply it to solve a problem in flight mechanics in a Newton's method context. In fact, a sequence of
linear problems is solved intending to  obtain a solution for the  original non-linear problem. We emphasize the method here
proposed has a performance considerably  better than those so far known, particularly concerning the computation time for a standard commercial  PC, which for the present one corresponds to a few seconds, even for a non-linear problem with about 40,000 degrees of freedom.

At this point we present a remark on the references.

\begin{remark} We highlight that a similar problem is addressed in \cite{12} for a nuclear physics model.
The main difference is that  now our results are more general and applicable to a much larger class of problems.
Specifically in the present work, we apply them to a flight mechanics model found in \cite{127}.

For the numerical results we have used finite differences. Details about finite differences schemes may be found in
\cite{103}.

Finally, details on the Sobolev spaces in which the original problem is established  may be found in \cite{1}.

\end{remark}

\section{The main results}

Consider the following system of difference equations in $\{(u_k)_n\}$, given by
\begin{equation}\label{a.1v}(u_k)_{n+1}=\sum_{j=1}^4 (a_{kj})_n(u_j)_n+(g_k)_n,\; \forall k \in \{1,2,3,4\},\; n \in \{0,...,N-1\},\end{equation}
where
$$(u_k)_n, (a_{kj})_n, (g_k)_n \in \mathbb{R},\; \forall j,k \in \{1,2,3,4\},\; n \in \{0,...,N-1\}.$$

Assume the following boundary conditions are intended to be satisfied:
\begin{equation}\left\{
\begin{array}{l}
 (u_1)_0=h_0,
 \\
 (u_3)_0=V_0 \\
 (u_4)_0=x_0 \\
 (u_1)_N=h_f. \end{array} \right. \end{equation}

For $n=0$ and $k=1$ we obtain
$$(u_1)_{1}=\sum_{j=1}^4( a_{1j})_0(u_j)_0+(g_1)_0,$$
that is
$$(u_2)_0=\frac{(u_1)_1-(g_1)_0-(a_{11})_0(u_1)_0-(a_{13})_0 (u_3)_0-(a_{14})_0 (u_4)_0 }{(a_{12})_0},$$
so that we write \begin{equation}\label{a.2v}(u_2)_0=m_2[0](u_1)_1 +\tilde{z}_2[0],\end{equation}
where,
$$m_2[0]=\frac{1}{(a_{12})_0},$$ and
$$\tilde{z}_2[0]=\frac{-(g_1)_0-(a_{11})_0(u_1)_0-(a_{13})_0 (u_3)_0-(a_{14})_0 (u_4)_0}{(a_{12})_0}.$$

Replacing (\ref{a.2v}) into
$$(u_2)_1=\sum_{j=1}^4 (a_{2j})_0 (u_j)_0+(g_2)_0,$$
we get
$$(u_2)_1=m_2[1](u_1)_1+z_2[1],$$
where
$$m_2[1]= (a_{22})_0 m_2[0],$$
and,
$$z_2[1]=(a_{21})_0(u_1)_0+(a_{22})_0 \tilde{z}_2[0]+(a_{23})_0 (u_3)_0+(a_{24})_0 (u_4)_0+(g_2)_0.$$

Also, replacing (\ref{a.2v}) into

$$(u_3)_1=\sum_{j=1}^4 (a_{3j})_0 (u_j)_0+(g_3)_0,$$ we obtain
$$(u_3)_1=m_3[1](u_1)_1+z_3[1],$$
where
$$m_3[1]= (a_{32})_0 m_2[0],$$
and,
$$z_3[1]=(a_{31})_0(u_1)_0+(a_{32})_0 \tilde{z}_2[0]+(a_{33})_0 (u_3)_0+(a_{34})_0 (u_4)_0+(g_3)_0.$$

Finally, replacing (\ref{a.2v}) into

$$(u_4)_1=\sum_{j=1}^4 (a_{4j})_0 (u_j)_0+(g_4)_0,$$
we may obtain
$$(u_4)_1=m_4[1](u_1)_1+z_4[1],$$
where
$$m_4[1]= (a_{42})_0 m_2[0],$$
and,
$$z_4[1]=(a_{41})_0(u_1)_0+(a_{42})_0 \tilde{z}_2[0]+(a_{43})_0 (u_3)_0+(a_{44})_0 (u_4)_0+(g_4)_0.$$

Reasoning inductively, having for $k \in \{2,3,4\},$
\begin{equation}\label{a.3v}(u_k)_{n-1}=m_k[n-1](u_1)_{n-1}+z_k[n-1],\end{equation}
for $n \geq 2$, replacing these last equations into (\ref{a.1v}) for $k=1$, we may obtain
\begin{equation}\label{a.4v}(u_1)_{n-1}=\tilde{m}_1[n-1](u_1)_n+\tilde{z}_1[n],\end{equation}
where
\begin{eqnarray}&&\tilde{m}_1[n-1]\nonumber \\ &=&\{(a_{11})_{n-1}+(a_{12})_{n-1}m_2[n-1]+(a_{13})_{n-1}m_3[n-1]
\nonumber \\ &&+(a_{14})_{n-1} m_4[n-1]\}^{-1},\end{eqnarray}
and
\begin{eqnarray}&&\tilde{z}_1[n-1] \nonumber \\ &=&-\tilde{m}_1[n-1]\{(a_{12})_{n-1}z_2[n-1]+(a_{13})_{n-1}z_3[n-1]\nonumber \\ &&+(a_{14})_{n-1}z_4[n-1]+(g_1)_{n-1}\}.\end{eqnarray}

Finally, replacing (\ref{a.4v}) into (\ref{a.3v}), we may obtain
\begin{equation}\label{a.5v}(u_k)_{n-1}=\tilde{m}_k[n-1](u_1)_n+\tilde{z}_k[n-1],\end{equation}
$\forall k \in \{2,3,4\},$

where
$$\tilde{m}_k[n-1]=m_k[n-1] \tilde{m}_1[n-1],$$
$$\tilde{z}_k[n-1]=m_k[n-1] \tilde{z}_1[n-1]+z_k[n-1].$$

Replacing (\ref{a.4v}) and  (\ref{a.5v}) into the system (\ref{a.1v}), we get
$$(u_k)_n=m_k[n](u_1)_n+z_k[n],$$
where,
$$m_k[n]=(a_{k1})_n \tilde{m}_1[n-1]+(a_{k2})_n \tilde{m}_2[n-1]+(a_{k3})_n \tilde{m}_3[n-1]+(a_{k4})_n \tilde{m}_4[n-1],$$
and,
$$z_k[n]=\sum_{j=1}^4 (a_{kj})_n\tilde{z}_j[n-1]+(g_k)_n.$$

Summarizing, we have obtained,
$$(u_1)_{n-1}=\tilde{m}_1[n-1](u_1)_n+\tilde{z}_1[n-1],$$
$$(u_k)_{n}=m_k[n](u_1)_n+z_k[n],\;\; \forall k \in \{2,3,4\}$$
$\forall n \in \{1,...,N\}.$

Therefore, having $(u_1)_N=h_f,$ we may obtain
$$(u_k)_N=m_k[N](u_1)_N+z_k[N],\;\; \forall k \in \{2,3,4\}$$
and
$$(u_1)_{N-1}=\tilde{m}_{N-1}(u_1)_N+\tilde{z}_1[N-1].$$

Having, $(u_1)_{N-1},$ we may obtain
$$(u_k)_{N-1}=m_k[N-1](u_1)_{N-1}+z_k[N-1],\;\; \forall k \in \{2,3,4\}$$
and
$$(u_1)_{N-2}=\tilde{m}_1[N-2](u_1)_{N-1}+\tilde{z}_1[N-2],$$
and so on, up to finding
$(u_k)_1, \forall k \in \{1,2,3,4\},$ and
finally,
$$(u_2)_0=m_2[0](u_1)_1+\tilde{z}_2[0].$$

At this point, the problem is solved.
\section{Numerical Results}

We present numerical results for the following system of equations, which models the in plane climbing  motion of
an airplane (please, see more details in \cite{127}).

\begin{equation}\left\{
\begin{array}{l}
 \dot{h}=V\sin \gamma,
 \\
 \dot{\gamma}=\frac{1}{m_f V}(T\sin(e_3)+L)-\frac{g}{V} \cos\gamma, \\
 \dot{V}=\frac{1}{m_f}(T\cos(e_3)-D)-g\sin\gamma \\
 \dot{x}=V\cos \gamma, \end{array} \right. \end{equation}

 with the boundary conditions,
 \begin{equation}\left\{
\begin{array}{l}
 h(0)=h_0,
 \\
 V(0)=V_0 \\
 x(0)=x_0 \\
 h(t_f)=h_f, \end{array} \right. \end{equation}

 where $t_f=40s$, $h$ is the airplane altitude, $V$ is its speed, $\gamma$ is the angle between its velocity and
 the horizontal axis, and finally $x$ denotes the horizontal coordinate position.

 For numerical purposes, we assume:

 $m_f=120,000 (2.2)\;lb,$ $S_f=260 (3.2)^2\;ft^2$, $a= 12^\circ$, $e_3=0.19$, $R_e=2.09\times 10^7\;ft$,
 $g=32\;ft/s^2$, $\rho_0=0.00239\; slug/ft^3,$ $B_9=26,600 \;ft,$ $\rho=\rho_0e^{-h/B_9},$
 $$C_L=-0.0005225 a^2+0.003506 a+0.1577(4.5),$$
 $$C_D=0.0001432 a^2+0.0058a+0.2204(1.8),$$
 $$L=\frac{1}{2}\rho V^2C_L S_f,$$
 $$D=\frac{1}{2}\rho V^2C_DS_f,$$
 and
 where units refer to the British system and,
 $$T=D+m_f g\sin\gamma,$$ which refers to a slightly negatively accelerated motion.

 To simplify the analysis, we redefine the variables as below indicated:

 \begin{equation}\left\{
\begin{array}{l}
 h=u_1,
 \\
 \gamma=u_2 \\
 V=u_3 \\
 x=u_4. \end{array} \right. \end{equation}

Thus, denoting $\mathbf{u}=(u_1,u_2,u_3,u_4) \in U=W^{1,2}([0,t_f];\mathbb{R}^4),$ the system above indicated may be expressed by

 \begin{equation}\label{a.6v}\left\{
\begin{array}{l}
 \dot{u}_1=f_1(\mathbf{u})
 \\
 \dot{u}_2=f_2(\mathbf{u})\\
 \dot{u}_3=f_3(\mathbf{u}) \\
 \dot{u}_4=f_4(\mathbf{u}), \end{array} \right. \end{equation}

where,

\begin{equation}\left\{
\begin{array}{l}
 f_1(\mathbf{u})=u_3\sin (u_2),
 \\
 f_2(\mathbf{u})=\frac{1}{m_f u_3}(T(\mathbf{u})\sin(e_3)+L(\mathbf{u}))-\frac{g}{u_3} \cos (u_2), \\
 f_3(\mathbf{u})=\frac{1}{m_f}(T(\mathbf{u})\cos(e_3)-D(\mathbf{u}))-g\sin (u_2) \\
f_4(\mathbf{u})= u_3 \cos (u_2). \end{array} \right. \end{equation}

Finally,
$$L(\mathbf{u})=\frac{1}{2}\rho_0e^{-u_1/B_9} u_3^2 C_L S_f,$$

$$D(\mathbf{u})=\frac{1}{2}\rho_0e^{-u_1/B_9} u_3^2 C_D S_f,$$

$$T(\mathbf{u})=D(\mathbf{u})+m_f g \sin (u_2).$$

At this point we shall write the system indicated in (\ref{a.6v}) in finite differences, that is,

 \begin{equation}\label{a.7v}\left\{
\begin{array}{l}
 (u_1)_{n+1}=(u_1)_n +f_1(\mathbf{u}_n)d
 \\
 (u_2)_{n+1}=(u_2)_n+f_2(\mathbf{u}_n)d\\
 (u_3)_{n+1}=(u_3)_n+f_3(\mathbf{u}_n)d \\
 (u_4)_{n+1}=(u_4)_n+f_4(\mathbf{u}_n)d, \end{array} \right. \end{equation}

 here $d=40/N,$ we $N$ refers to the number of nodes concerning the discretization in $t$ (in our numerical
 example $N=10000$).

 Intending to apply the Newton's method we linearize the system indicated in (\ref{a.7v}) about
 a initial guess $$\tilde{\mathbf{u}}=(\tilde{u}_1,\tilde{u}_2,\tilde{u}_3, \tilde{u}_4).$$

 We obtain the following approximate system
 \begin{equation}\label{a.8v}\left\{
\begin{array}{l}
 (u_1)_{n+1}=(u_1)_n +d\;\left(f_1(\tilde{\mathbf{u}}_n)+\sum_{j=1}^4 \frac{\partial f_1(\tilde{\mathbf{u}}_n)}{\partial \tilde{u}_j}((u_j)_n-(\tilde{u}_j)_n)\right)
 \\
 (u_2)_{n+1}=(u_2)_n+d\;\left(f_2(\tilde{\mathbf{u}}_n)+\sum_{j=1}^4 \frac{\partial f_2(\tilde{\mathbf{u}}_n)}{\partial \tilde{u}_j}((u_j)_n-(\tilde{u}_j)_n)\right)\\
 (u_3)_{n+1}=(u_3)_n+d\;\left(f_3(\tilde{\mathbf{u}}_n)+\sum_{j=1}^4 \frac{\partial f_3(\tilde{\mathbf{u}}_n)}{\partial \tilde{u}_j}((u_j)_n-(\tilde{u}_j)_n)\right)\\
 (u_4)_{n+1}=(u_4)_n+d\;\left(f_4(\tilde{\mathbf{u}}_n)+\sum_{j=1}^4 \frac{\partial f_4(\tilde{\mathbf{u}}_n)}{\partial \tilde{u}_j}((u_j)_n-(\tilde{u}_j)_n)\right). \end{array} \right. \end{equation}

Observe that such a system is in the form,

$$(u_k)_{n+1}=(u_k)_n+\sum_{j=1}^4 (a_{kj})_n(u_j)_n+(g_k)_n,$$
where
$$(a_{kj})_n=\frac{\partial f_k(\tilde{\mathbf{u}}_n)}{\partial u_j} \;d, \text{ for } j \neq k$$

$$(a_{jj})_n=1+\frac{\partial f_j(\tilde{\mathbf{u}}_n)}{\partial u_j} \;d, \text{ for } j=k,$$

and,
$$(g_k)_n=f_k(\tilde{\mathbf{u}}_n)\;d-\sum_{j=1}^4\frac{\partial f_k(\tilde{\mathbf{u}}_n)}{\partial u_j}(\tilde{u}_j)_n \;d.$$

We solve this last system for the following boundary conditions:
 \begin{equation}\left\{
\begin{array}{l}
 h(0)=0\; ft,
 \\
 V(0)=960\; ft/s, \\
 x(0)=0\;ft, \\
 h(t_f)=35000\;ft. \end{array} \right. \end{equation}

We have obtained $\{\mathbf{u}_n\}$. In a Newton's method context, the next step is to replace
$\tilde{\mathbf{u}}_n$ by   $\{\mathbf{u}_n\}$ and  thus to repeat the process up to the satisfaction of an
appropriate convergence criterion.

\begin{figure}
\centering \includegraphics [width=3in]{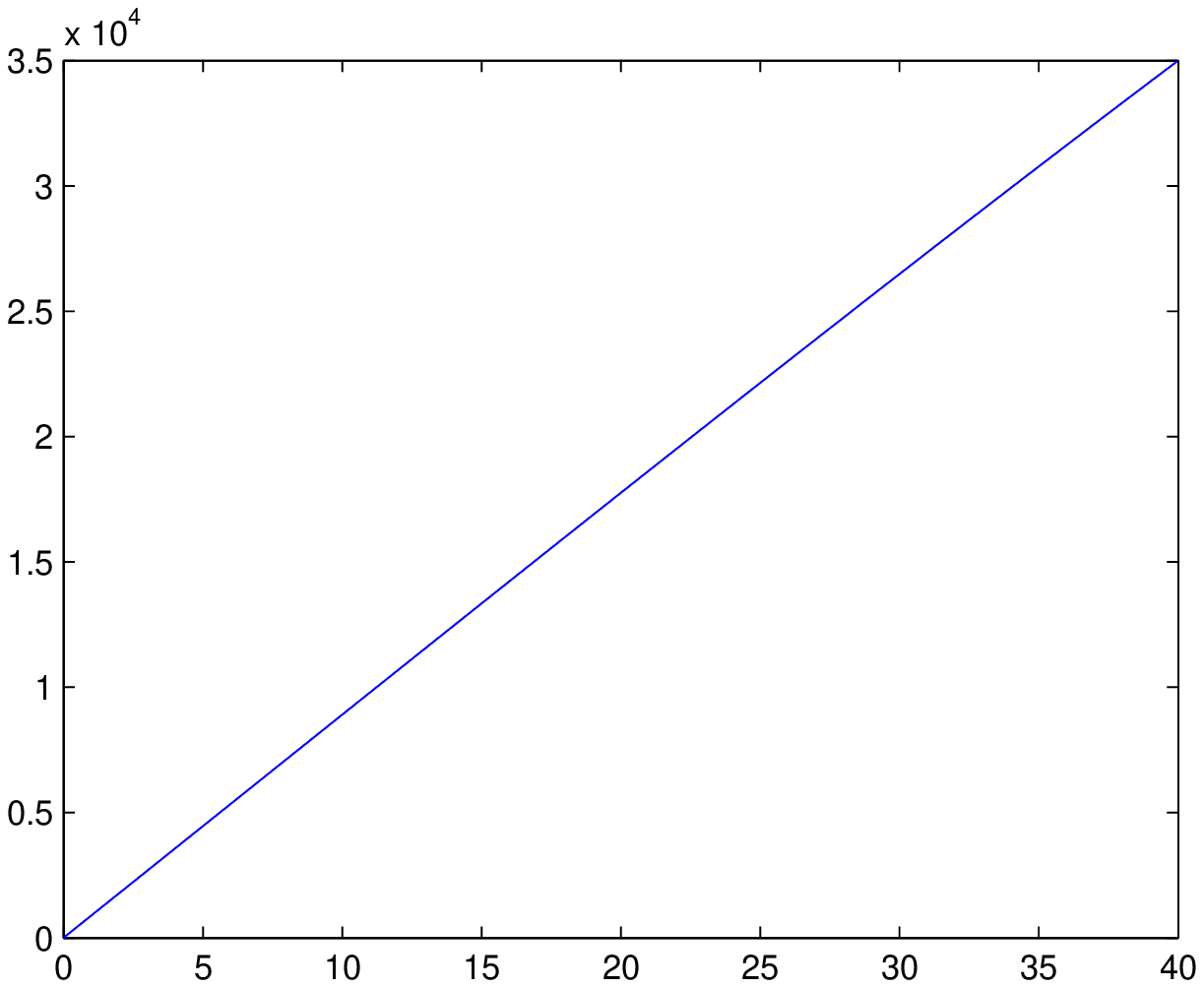}
\\ \caption{\small{The solution $h$ (in ft) for $t_f=40s$.
}}\label{fig.1v}
\end{figure}

\begin{figure}
\centering \includegraphics [width=3in]{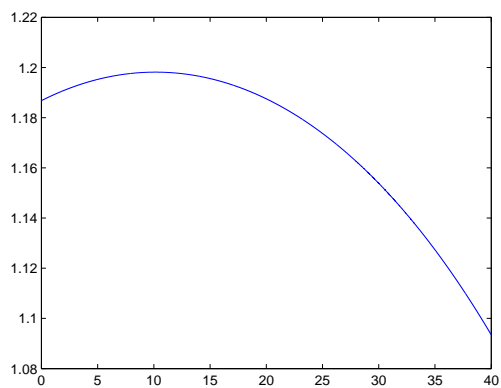}
\\ \caption{\small{The solution $\gamma$ (in rad) for $t_f=40s$.
}}\label{fig.2v}
\end{figure}

\begin{figure}
\centering \includegraphics [width=3in]{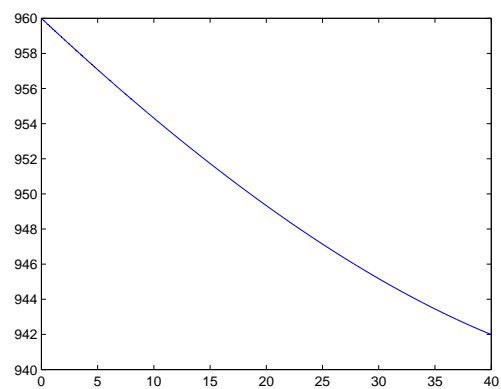}
\\ \caption{\small{The solution $V$ (in ft/s) for $t_f=40s$.
}}\label{fig.3v}
\end{figure}

\begin{figure}
\centering \includegraphics [width=3in]{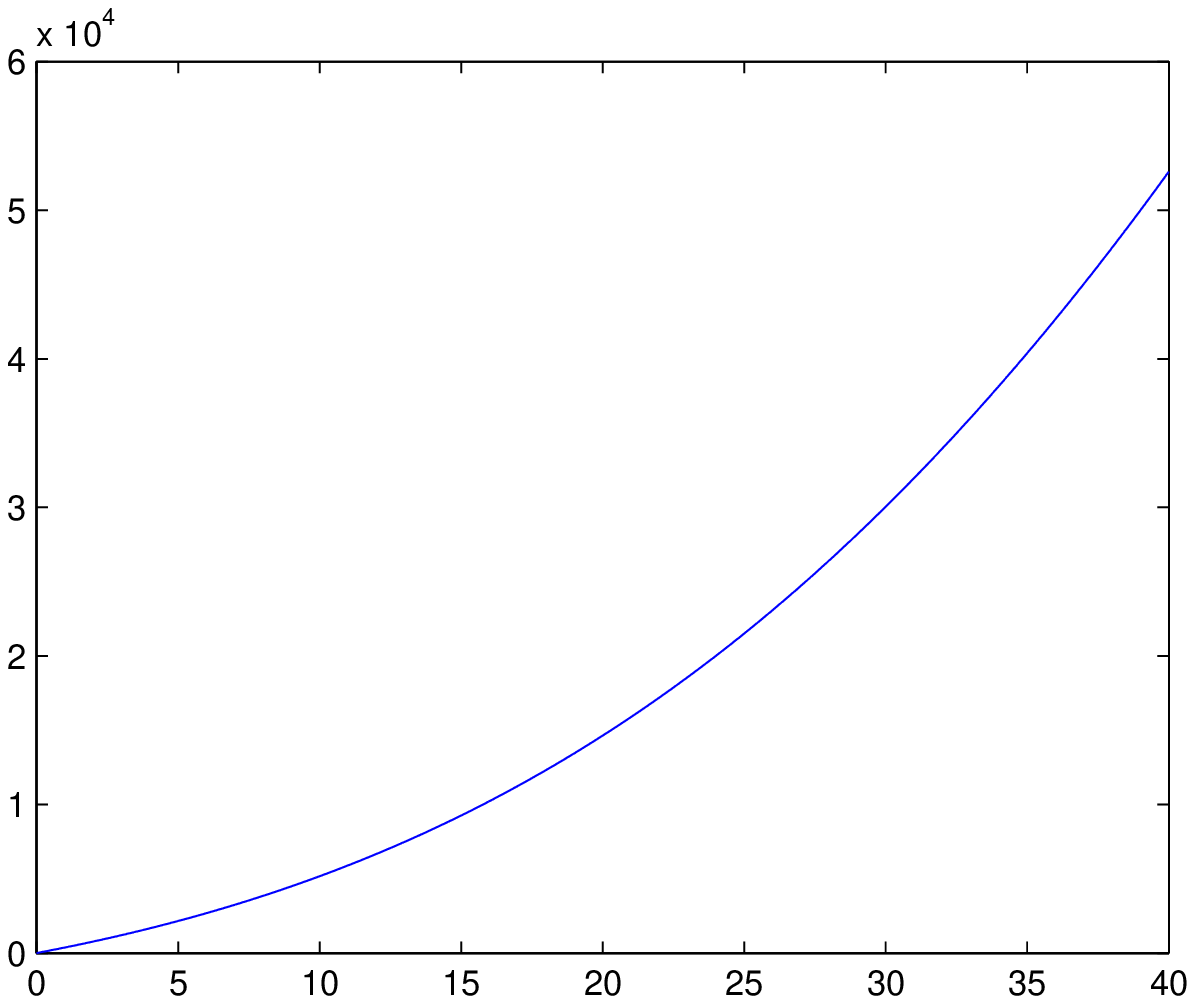}
\\ \caption{\small{The solution $x$ (in ft) for $t_f=40s$.
}}\label{fig.4v}
\end{figure}

We have obtained the following solutions for $h,\gamma,V \text{ and } x$. Please see figures \ref{fig.1v}, \ref{fig.2v},  \ref{fig.3v} and \ref{fig.4v}, respectively.


\section{Conclusion}
In this article, we have developed a method for solving a class of first order
ordinary differential equations.

The results are applied to a flight mechanics problem
which models the in plane climbing of an airplane. It is worth mentioning the algorithm obtained is of
easy implementation and very efficient from a computational point of view.

Finally, we would highlight the
numerical results obtained are perfectly consistent with the physical problem context.
In future works we intend to apply
the method to  solve relating optimal control problems.

%Non-BibTeX users please useu


\begin{thebibliography}{}
\bibitem{1} R.A. Adams and J.F. Fournier, Sobolev Spaces, 2nd edn. Elsevier, New York, 2003.

\bibitem{12} F. Botelho, Functional Analysis and Applied Optimization in Banach Spaces,
 Springer Switzerland, 2014.

\bibitem{103}  J.C. Stikwerda, Finite Difference Schemes and Partial Differential Equations, 2nd edn. SIAM, Philadephia, 2004.

\bibitem{127} N.X. Vinh, Flight Mechanics of High Performance Aircraft, Cambridge University Press, New York, 1993.
\end{thebibliography}
\end{document}